\documentclass[12pt,a4paper,oneside,reqno]{article}

\usepackage[utf8]{inputenc}
\usepackage[english]{babel}
\usepackage[left=1in,right=1in,top=1in,bottom=1in]{geometry}
\usepackage{amsfonts,amssymb,amsmath,amsthm}
\usepackage{enumerate}
\usepackage{cite}
\usepackage{xcolor}
\usepackage[normalem]{ulem}
\usepackage{longtable}
\usepackage{hyperref}
\hypersetup{colorlinks=true,linkcolor=blue,citecolor=blue,urlcolor=blue}
\usepackage{fancyhdr}
\theoremstyle{plain}
\newtheorem{Th}{Theorem} 
\newtheorem{Lem}[Th]{Lemma}
\theoremstyle{remark}
\newtheorem*{RemA}{Remark A (complete multiplication table of $K_{10}$)}

\newcommand{\FF}{F}
\newcommand{\ov}[1]{\overline{#1}}

\newcommand{\doi}[1]{\href{https://doi.org/#1}{doi:#1}}

\begin{document}
\noindent\fbox{\parbox{\dimexpr\textwidth-2\fboxsep-2\fboxrule\relax}{\small
\textbf{Corrected version verified in Lean.} This is a corrected version of the paper
I.~Kaygorodov, \emph{$\delta$-derivations of simple finite-dimensional Jordan superalgebras},
Algebra and Logic 46 (2007), no.~5, 318--329, \doi{10.1007/s10469-007-0032-0}.}}

\bigskip

\begin{center}
{\Large $\delta$-derivations  of 
simple finite-dimensional Jordan superalgebras}

\bigskip

{\bf Ivan~Kaygorodov}
\end{center}

\bigskip

\noindent {\bf Abstract.} {\it We describe non-trivial $\delta$-derivations of semisimple finite-dimensional
Jordan algebras over an algebraically closed field of characteristic not $2$, and of simple
finite-dimensional Jordan superalgebras over an algebraically closed field of characteristic $0$. For
these classes of algebras and superalgebras, non-zero non-trivial $\delta$-derivations are shown
to be missing for $\delta\notin\{ 0,\frac12,1\}$, and we give a complete account of $\frac12$-derivations.}

\bigskip


\bigskip
\setcounter{tocdepth}{1}
\tableofcontents

\bigskip

\section*{Introduction}
\addcontentsline{toc}{section}{Introduction}

The notion of derivation for an algebra was generalized by many mathematicians along quite different
lines. Thus, in \cite{Jac62}, the reader can find the definitions of a derivation of a subalgebra into an
algebra and of an $(s_1,s_2)$-derivation of one algebra into another, where $s_1$ and $s_2$ are some
homomorphisms of the algebras. Back in the 1950s, Herstein explored Jordan derivations of prime
associative rings of characteristic $p\neq 2$; see \cite{Her57}. (Recall that a Jordan derivation of an
algebra $A$ is a linear mapping $j_d\colon A\to A$ satisfying the equality
$j_d(xy+yx)=j_d(x)y+xj_d(y)+j_d(y)x+yj_d(x)$ for any $x,y\in A$.) He proved that the Jordan derivation
of such a ring is properly a standard derivation. Later on, Hopkins in \cite{Hop96} dealt with
antiderivations of Lie algebras (for the definition of an antiderivation, see \cite{Jac62}). The
antiderivation, on the other hand, is a special case of a $\delta$-derivation --- that is, a linear
mapping $\mu$ of an algebra such that $\mu(xy)=\delta(\mu(x)y+x\mu(y))$, where $\delta$ is some fixed
element of the ground field.

Subsequently, Filippov generalized Hopkins' results in \cite{Fil98} and
\cite{Fil99} by treating prime Lie algebras over an associative commutative ring $\Phi$ with unity and
$\frac12$. It was proved that every prime Lie $\Phi$-algebra, on which a non-degenerate symmetric
invariant bilinear form is defined, has no non-zero $\delta$-derivation if $\delta\neq -1,0,\frac12,1$.
In \cite{Fil99}, also, $\frac12$-derivations were described for an arbitrary prime Lie $\Phi$-algebra $A$
($\frac16\in\Phi$) with a non-degenerate symmetric invariant bilinear form defined on the algebra. It
was shown that the linear mapping $\varphi\colon A\to A$ is a $\frac12$-derivation iff
$\varphi\in\Gamma(A)$, where $\Gamma(A)$ is the centroid of $A$. This implies that if $A$ is a central
simple Lie algebra over a field of characteristic $p\neq 2,3$ on which a non-degenerate symmetric
invariant bilinear form is defined, then every $\frac12$-derivation $\varphi$ has the form
$\varphi(x)=\alpha x$, $\alpha\in\Phi$. At a later time, Filippov described $\delta$-derivations for prime
alternative and non-Lie Malcev $\Phi$-algebras with some restrictions on the operator
ring $\Phi$. In \cite{Fil00}, for instance, it was stated that algebras in these classes have no non-zero
$\delta$-derivations if $\delta\neq 0,\frac12,1$.

In the present paper, we come up with an account of non-trivial $\delta$-derivations for semisimple
finite-dimensional Jordan algebras over an algebraically closed field of characteristic not $2$, and for
simple finite-dimensional Jordan superalgebras over an algebraically closed field of characteristic $0$.
For these classes of algebras and superalgebras, non-zero non-trivial $\delta$-derivations are
shown to be missing for $\delta\neq 0,\frac12,1$, and we provide  a complete description of
$\frac12$-derivations.

The paper is divided into four parts. In Sec.~1, relevant definitions are given and known results cited.
In Sec.~2, we deal with $\delta$-derivations of simple and semisimple finite-dimensional Jordan algebras.
In Secs.~3 and~4, $\delta$-derivations are described for simple finite-dimensional Jordan
superalgebras over an algebraically closed field of characteristic $0$. For
some superalgebras, note, the condition on the characteristic may be weakened so as to be distinct from
$2$. A proof for the main theorem is based on the classification theorem for simple finite-dimensional
Jordan superalgebras and on the results obtained in Secs.~3 and~4.

\section{Basic facts and definitions}

Let $\FF$ be a field of characteristic $p$, $p\neq 2$. An algebra $A$ over $\FF$ is Jordan if it
satisfies the following identities:
\[xy=yx,\qquad (x^2y)x=x^2(yx).\]
Jordan algebras arise naturally from the associative algebras. If in an associative algebra $A$ we
replace multiplication $ab$ by symmetrized multiplication $a\circ b=\frac12(ab+ba)$ then we will face a
Jordan algebra. Denote this algebra by $A^{(+)}$. Below are essential examples of Jordan algebras.

(1) The algebra $J(V,f)$ of a bilinear form. Let $f\colon V\times V\to \FF$ be a symmetric bilinear form
on a vector space $V$. On the direct sum $J=\FF\cdot 1+V$ of vector spaces, we then define
multiplication by setting $1\cdot v=v\cdot 1=v$ and $v_1\cdot v_2=f(v_1,v_2)\cdot 1$; under this
multiplication, $J=J(V,f)$ is a Jordan algebra. If the form $f$ is non-degenerate and $\dim V>1$, then
the algebra $J(V,f)$ is simple.

(2) The Jordan algebra $H(D_n,J)$. Here, $n\geqslant 3$, $D$ is a composition algebra, which is
associative for $n>3$, $j\colon d\mapsto\bar d$ is a canonical involution in $D$, and $J\colon X\mapsto
\bar X^{T}$ is a standard involution in $D_n$.

\begin{Th}[see \cite{ZSSS78}]\label{th11}
Every simple finite-dimensional Jordan algebra $A$ over an algebraically closed field $\FF$ of
characteristic not $2$ is isomorphic to one of the following algebras:
(1) $\FF\cdot 1$; (2) $J(V,f)$; (3) $H(D_n,J)$.
\end{Th}

We recall the definition of a superalgebra. Let $\Gamma$ be a Grassmann algebra over $\FF$, which is
generated by elements $1,e_1,\dots,e_n,\dots$ and is defined by relations $e_i^2=0$,
$e_ie_j=-e_je_i$. Products $1,e_{i_1}e_{i_2}\dots e_{i_k}$, $i_1<i_2<\dots<i_k$, form a basis for
$\Gamma$ over $\FF$. Denote by $\Gamma_0$ and $\Gamma_1$ the subspaces generated by products of even
and odd lengths, respectively. Then $\Gamma$ is represented as a direct sum of these subspaces,
$\Gamma=\Gamma_0+\Gamma_1$, with $\Gamma_i\Gamma_j\subseteq\Gamma_{i+j\,(\mathrm{mod}\,2)}$,
$i,j=0,1$. In other words, $\Gamma$ is a $\mathbb Z_2$-graded algebra (or superalgebra) over $\FF$.

Now let $A=A_0+A_1$ be any superalgebra over $\FF$. Consider a tensor product
of $\FF$-algebras, $\Gamma\otimes A$. Its subalgebra
\[\Gamma(A)=\Gamma_0\otimes A_0+\Gamma_1\otimes A_1\]
is called a Grassmann envelope for $A$.

Let $\Omega$ be some variety of algebras over $\FF$. A $\mathbb Z_2$-graded algebra $A=A_0+A_1$ is an
$\Omega$-superalgebra if its Grassmann envelope $\Gamma(A)$ is an algebra in $\Omega$. In particular,
$A=A_0\oplus A_1$ is a Jordan superalgebra if its Grassmann envelope $\Gamma(A)$ is a Jordan algebra.

In \cite{Wal64}, it was shown that every simple finite-dimensional associative superalgebra over an
algebraically closed field $\FF$ is isomorphic either to $A=M_{m,n}(\FF)$, which is the matrix algebra
$M_{m+n}(\FF)$, or to $B=Q(n)$, which is a subalgebra of $M_{2n}(\FF)$. Gradings of the superalgebras $A$
and $B$ are the following:
\begin{longtable}{rcl}
$A_0$&$=$&$\left\{\begin{pmatrix}a&0\\0&d\end{pmatrix}\ \Big|\ a\in M_m(\FF),\ d\in M_n(\FF)\right\},$\\
$A_1$&$=$&$\left\{\begin{pmatrix}0&b\\c&0\end{pmatrix}\ \Big|\ b\in M_{m,n}(\FF),\ c\in M_{n,m}(\FF)\right\},$\\
$B_0$&$=$&$\left\{\begin{pmatrix}a&0\\0&a\end{pmatrix}\ \Big|\ a\in M_n(\FF)\right\},$\\
$B_1$&$=$&$\left\{\begin{pmatrix}0&b\\b&0\end{pmatrix}\ \Big|\ b\in M_n(\FF)\right\}.$
\end{longtable}

Let $A=A_0+A_1$ be an associative superalgebra. The vector space of $A$ can be endowed with the
structure of a Jordan superalgebra $A^{(+)}$, by defining new multiplication as
follows: $a\circ b=\frac12(ab+(-1)^{p(a)p(b)}ba)$. In this case $p(a)=i$ if $a\in A_i$.

Using the above construction, we arrive at superalgebras
\[M_{m,n}(\FF)^{(+)},\ m\geqslant 1,\ n\geqslant 1;\qquad Q(n)^{(+)},\ n\geqslant 2.\]

Now, we define the superinvolution $j\colon A\to A$. A graded endomorphism $j\colon A\to A$ is called a
superinvolution if $j(j(a))=a$ and $j(ab)=(-1)^{p(a)p(b)}j(b)j(a)$. Let $H(A,j)=\{a\in A: j(a)=a\}$.
Then $H(A,j)=H(A_0,j)+H(A_1,j)$ is a subsuperalgebra of $A^{(+)}$. Below are superalgebras which are
obtained from $M_{n,m}(\FF)$ via a suitable superinvolution:

(1) the Jordan superalgebra $osp(n,m)$, consisting of matrices of the form
$\begin{pmatrix}A&B\\C&D\end{pmatrix}$, where $A^T=A\in M_n(\FF)$, $C=Q^{-1}B^T$,
$D=Q^{-1}D^TQ\in M_{2m}(\FF)$, and $Q=\begin{pmatrix}0&E_m\\-E_m&0\end{pmatrix}$;

(2) the Jordan superalgebra $P(n)$, consisting of matrices of the form
$\begin{pmatrix}A&B\\C&D\end{pmatrix}$, where $B^T=-B$, $C^T=C$, and $D=A^T$, with
$A,B,C,D\in M_n(\FF)$.

\begin{Th}[see  \cite{Kan89,Kac77,RZ94}]\label{th12}
Every simple finite-dimensional non-trivial (i.e., with a non-zero odd part) Jordan superalgebra $A$ over
an algebraically closed field $\FF$ of characteristic $0$ is isomorphic to one of the following
superalgebras:
\[M_{m,n}(\FF)^{(+)};\ Q(n)^{(+)};\ osp(n,m);\ P(n);\ J(V,f);\ D_t,\ t\neq 0;\ K_3;\ K_{10};\ J(\Gamma_n),\
n>1.\]
\end{Th}

The superalgebras $J(V,f)$, $D_t$, $K_3$, $K_{10}$, and $J(\Gamma_n)$ will be defined below.

Let $\delta\in\FF$. A linear mapping $\varphi$ of $A$ is called a $\delta$-derivation if
\begin{equation}\label{eq1}
\varphi(xy)=\delta(x\varphi(y)+\varphi(x)y)
\end{equation}
for arbitrary elements $x,y\in A$.

The definition of a $1$-derivation coincides with the conventional definition of a derivation. A
$0$-derivation is any endomorphism $\varphi$ of $A$ such that $\varphi(A^2)=0$. A non-trivial
$\delta$-derivation is a $\delta$-derivation which is not a $1$-derivation, nor a $0$-derivation.
Obviously, for any algebra, the multiplication operator by an element of the ground field $\FF$ is a
$\frac12$-derivation. We are interested in the behavior of non-trivial $\delta$-derivations of semisimple
finite-dimensional Jordan algebras over an algebraically closed field of characteristic not $2$, and of
simple finite-dimensional Jordan superalgebras over an algebraically closed field of characteristic $0$.

\section{$\delta$-Derivations for semisimple finite-dimensional Jordan algebras}

In this section, we look at how non-trivial $\delta$-derivations of simple finite-dimensional Jordan
algebras behave over an algebraically closed field $\FF$ of characteristic distinct from $2$. As a
consequence, we furnish a description of $\delta$-derivations for semisimple finite-dimensional Jordan
algebras over an algebraically closed field of characteristic not $2$.

\begin{Th}\label{th21}
Let $\varphi$ be a non-trivial $\delta$-derivation of a superalgebra $A$ with unity $e$ over a field $\FF$
of characteristic not $2$. Then $\delta=\frac12$.
\end{Th}

\begin{proof}
Let $\delta\neq\frac12$. Then $\varphi(e)=\varphi(e\cdot e)=\delta(\varphi(e)+\varphi(e))=2\delta\varphi(e)$,
that is, $\varphi(e)=0$. Thus $\varphi(x)=\varphi(x\cdot e)=\delta(\varphi(x)+x\varphi(e))=\delta\varphi(x)$
for arbitrary $x\in A$. Hence either $\delta=1$, or $\varphi=0$, and in both cases $\varphi$ is
trivial. Contradiction. The theorem is proved.
\end{proof}

\begin{Lem}\label{lem22}
Let $\varphi$ be a non-trivial $\frac12$-derivation of a Jordan algebra $A$ isomorphic to the ground field.
Then $\varphi(x)=\alpha x$, $\alpha\in\FF$.
\end{Lem}

\begin{proof}
Let $e$ be unity in $A$. Then
\begin{equation}\label{eq2}
\varphi(x)=2\varphi(xe)-\varphi(x)=x\varphi(e),
\end{equation}
that is, $\varphi(x)=\alpha x$, $\alpha\in\FF$. The lemma is proved.
\end{proof}

Note that identity~\eqref{eq2} holds for every $\frac12$-derivation of every (super)algebra with
unity $e$; it is used in this form below.

\begin{Lem}\label{lem23}
Let $\varphi$ be a non-trivial $\frac12$-derivation of an algebra $J(V,f)$, where $f$ is
non-degenerate and $\dim V>1$. Then $\varphi(x)=\alpha x$ for $\alpha\in\FF$.
\end{Lem}

\begin{proof}
Let $\varphi(e)=\alpha e+v$, where $\alpha\in\FF$ and $v\in V$. From \eqref{eq2}, it follows that
$\varphi(x)=x\varphi(e)$ for any $x\in J(V,f)$.

For $w\in V$, we then have
\begin{multline*}
\alpha f(w,w)e+f(w,w)v=w^2(\alpha e+v)=\varphi(w^2)=\tfrac12(w\varphi(w)+\varphi(w)w)\\
=w\varphi(w)=w(w(\alpha e+v))=w(\alpha w+f(v,w)e)=\alpha f(w,w)e+f(w,v)w.
\end{multline*}
As the result, $f(w,w)v=f(w,v)w$. Now, since $w$ is arbitrary and $\dim(V)>1$, we have $v=0$. Thus
$\varphi(x)=\alpha x$ for any $x\in J(V,f)$. The lemma is proved.
\end{proof}

\begin{Lem}\label{lem24}
Let $\varphi$ be a non-trivial $\frac12$-derivation of an algebra $H(D_n,J)$, $n\geqslant 3$. Then
$\varphi(x)=\alpha x$ for $\alpha\in\FF$.
\end{Lem}

\begin{proof}
Relevant information on composition algebras can be found in \cite{ZSSS78}. Let $\varphi(e)=\alpha e+v$,
where $v=\sum\limits_{i,j=1}^{n}x_{i,j}e_{i,j}$, $x_{1,1}=0$,
$x_{j,i}=\ov{x_{i,j}}$, $\alpha\in\FF$, $x_{i,j}\in D$
(the diagonal entries $x_{i,i}$ are symmetric elements of $D$, i.e.\ elements of $\FF$).

From \eqref{eq2}, for $x\in H(D_n,J)$ arbitrary, we have
\begin{equation}\label{eq3}
x^2\circ(\alpha e+v)=\varphi(x^2)=x\circ\varphi(x)=x\circ(x\circ(\alpha e+v)),\qquad x^2\circ v=x\circ(x\circ v).
\end{equation}
If we put $x=e_{k,k}$ we obtain
\begin{center}
    $\sum_{j=1}^n x_{k,j}e_{k,j}+\sum_{i=1}^n x_{i,k}e_{i,k}=2e_{k,k}^2\circ v=2e_{k,k}\circ(e_{k,k}\circ v)
=$\\
$\frac12\Big(\sum_{j=1}^n x_{k,j}e_{k,j}+x_{k,k}e_{k,k}+x_{k,k}e_{k,k}+\sum_{i=1}^n x_{i,k}e_{i,k}\Big),$
\end{center}
whence $v=\sum\limits_{i=1}^n x_{i,i}e_{i,i}$.

For $x=e_{n,k}+e_{k,n}$ substituted in \eqref{eq3}, we have
{\small \begin{multline*}
x_{n,n}e_{n,n}+x_{k,k}e_{k,k}=(e_{n,k}+e_{k,n})^2\circ\sum_{i=1}^n x_{i,i}e_{i,i}
=(e_{n,k}+e_{k,n})\circ\Big((e_{n,k}+e_{k,n})\circ\sum_{i=1}^n x_{i,i}e_{i,i}\Big)
=\\(e_{n,k}+e_{k,n})\circ\tfrac12(x_{n,n}e_{k,n}+x_{k,k}e_{k,n}+x_{k,k}e_{n,k}+x_{n,n}e_{n,k})
=\\\tfrac12(x_{k,k}e_{k,k}+x_{k,k}e_{n,n}+x_{n,n}e_{k,k}+x_{n,n}e_{n,n}),
\end{multline*}}
which yields $x_{n,n}=x_{n-1,n-1}=\dots=x_{1,1}=0$ and $v=0$.

Consequently, $\varphi(x)=\alpha x$ for any $x\in H(D_n,J)$. The lemma is proved.
\end{proof}

\begin{Th}\label{th25}
Let $\varphi$ be a non-trivial $\delta$-derivation of a simple finite-dimensional Jordan algebra $A$ over
an algebraically closed field $\FF$ of characteristic distinct from $2$. Then $\delta=\frac12$ and
$\varphi(x)=\alpha x$, $\alpha\in\FF$.
\end{Th}

The proof follows from Theorems~\ref{th11}, \ref{th21} and Lemmas~\ref{lem22}--\ref{lem24}.

\begin{Th}\label{th26}
Let $\varphi$ be a non-trivial $\delta$-derivation of a semisimple finite-dimensional Jordan algebra
$A=\bigoplus\limits_{i=1}^n A_i$, where $A_i$ are simple algebras, over an algebraically closed field of
characteristic not $2$. Then $\delta=\frac12$, and for $x=\sum\limits_{i=1}^n x_i$ where $x_i\in A_i$, we
have $\varphi(x)=\sum\limits_{i=1}^n\alpha_ix_i$, $\alpha_i\in\FF$.
\end{Th}

\begin{proof}
The algebra $A$ has the unity $\sum_k e_k$, so $\delta=\frac12$ by Theorem~\ref{th21}. Unity in
$A_k$ is denoted by $e_k$. If $x_i\in A_i$, then $\varphi(x_i)=x_i^++x_i^-$, where $x_i^+\in A_i$ and
$x_i^-\in\sum_{k\neq i}A_k$. Put $e^i=\sum\limits_{k=1}^n e_k-e_i$ and $\varphi(e^i)=e^{i+}+e^{i-}$, where
$e^{i+}\in A_i$ and $e^{i-}\in\sum_{k\neq i}A_k$. Then
\[0=\varphi(x_i\cdot e^i)=\delta(\varphi(x_i)\cdot e^i+x_i\cdot\varphi(e^i))
=\delta((x_i^++x_i^-)e^i+x_i(e^{i+}+e^{i-}))=\delta(x_i^-+x_i\cdot e^{i+}),\]
which yields $x_i^-=0$. Consequently, the mapping $\varphi$ is invariant on $A_i$. In virtue of
Theorem~\ref{th25},  $\varphi(x_i)=\alpha_ix_i$ for some $\alpha_i\in\FF$ defined
for $A_i$ with $x_i\in A_i$ arbitrary (the restriction of $\varphi$ to $A_i$ is a
$\frac12$-derivation of $A_i$, and every $\frac12$-derivation of $A_i$ is scalar by Theorem~\ref{th25}
and the remark preceding Theorem~\ref{th21}). It is easy to verify that the mapping $\varphi$, given by
the rule $\varphi\big(\sum_{i=1}^n x_i\big)=\sum_{i=1}^n\alpha_ix_i$, $x_i\in A_i$, is a
$\frac12$-derivation. The theorem is proved.
\end{proof}

\section{$\delta$-Derivations for simple finite-dimensional Jordan superalgebras with unity}

In this section, all superalgebras but $J(\Gamma_n)$ are treated over a field of characteristic not $2$.
The superalgebra $J(\Gamma_n)$ is treated over a field of characteristic $0$. Among the title
superalgebras are $M_{m,n}(\FF)^{(+)}$, $Q(n)^{(+)}$, $osp(n,m)$, $P(n)$, $J(V,f)$, and $J(\Gamma_n)$.
Theorem~\ref{th21} implies that these superalgebras all lack  non-trivial $\delta$-derivations,
for $\delta\neq\frac12$. Therefore, we need only consider the case of a $\frac12$-derivation.

\begin{Lem}\label{lem31}
Let $\varphi$ be a non-trivial $\frac12$-derivation of $M_{m,n}(\FF)^{(+)}$. Then $\varphi(x)=\alpha x$
for some $\alpha\in\FF$.
\end{Lem}

\begin{proof}
It is easy to see that, for $1\leqslant i,j\leqslant n+m$, elements $e_{i,j}$ form a basis for the
superalgebra $M_{m,n}(\FF)^{(+)}$. Let $\varphi(e_{i,j})=\sum\limits_{k,l=1}^{m+n}\alpha^{i,j}_{k,l}e_{k,l}$,
where $\alpha^{i,j}_{k,l}\in\FF$, $i,j=1,\dots,n+m$.

If in \eqref{eq1} we put $x=y=e_{i,i}$ we arrive at
\[\sum_{k,l=1}^{m+n}\alpha^{i,i}_{k,l}e_{k,l}=\varphi(e_{i,i})=\varphi(e_{i,i}^2)
=\tfrac12(e_{i,i}\circ\varphi(e_{i,i})+\varphi(e_{i,i})\circ e_{i,i})
=\tfrac12\Big(\sum_{l=1}^{n+m}\alpha^{i,i}_{i,l}e_{i,l}+\sum_{k=1}^{n+m}\alpha^{i,i}_{k,i}e_{k,i}\Big),\]
whence $\varphi(e_{i,i})=\alpha_ie_{i,i}$, where $\alpha_i=\alpha^{i,i}_{i,i}$, $i=1,\dots,m+n$.

Substituting $x=e_{i,j}$ and $y=e_{i,i}$, $i\neq j$, in \eqref{eq1}, we obtain
\[\sum_{k,l=1}^{m+n}\alpha^{i,j}_{k,l}e_{k,l}=\varphi(e_{i,j})=2\varphi(e_{i,j}\circ e_{i,i})
=\tfrac12\Big(\alpha_ie_{i,j}+\sum_{l=1}^{m+n}\alpha^{i,j}_{i,l}e_{i,l}+\sum_{k=1}^{m+n}\alpha^{i,j}_{k,i}e_{k,i}\Big).\]
Analyzing the resulting equalities, we conclude that $\alpha^{i,j}_{i,j}=\alpha_i$. A similar argument for
$e_{i,j}$ and $e_{j,j}$ yields $\alpha^{i,j}_{i,j}=\alpha_j$. Since $\varphi$ is linear, $\varphi(e)=\alpha e$.
Using \eqref{eq2} gives $\varphi(x)=\alpha x$, for any $x\in M_{m,n}(\FF)^{(+)}$. The lemma
is proved.
\end{proof}

\begin{Lem}\label{lem32}
Let $\varphi$ be a non-trivial $\frac12$-derivation of $Q(n)^{(+)}$. Then $\varphi(x)=\alpha x$, where
$\alpha\in\FF$.
\end{Lem}

\begin{proof}
Clearly, $\Delta_{i,j}=e_{i,j}+e_{n+i,n+j}$ and $\Delta^{i,j}=e_{n+i,j}+e_{i,n+j}$ form a basis for the
superalgebra $Q(n)^{(+)}$.

On the basis elements, the following relations hold:
\[\Delta_{i,j}\circ\Delta_{k,l}=\tfrac12(\delta_{j,k}\Delta_{i,l}+\delta_{l,i}\Delta_{k,j}),\qquad
\Delta_{i,j}\circ\Delta^{k,l}=\tfrac12(\delta_{j,k}\Delta^{i,l}+\delta_{l,i}\Delta^{k,j}).\]
Let $\varphi(\Delta_{i,j})=\sum\limits_{k,l=1}^n\alpha^{i,j}_{k,l}\Delta_{k,l}
+\sum\limits_{k,l=1}^n\alpha^{*i,j}_{k,l}\Delta^{k,l}$. Put $x=y=\Delta_{i,i}$ in \eqref{eq1}. Then
\begin{multline*}
\sum_{k,l=1}^n\alpha^{i,i}_{k,l}\Delta_{k,l}+\sum_{k,l=1}^n\alpha^{*i,i}_{k,l}\Delta^{k,l}
=\varphi(\Delta_{i,i})=\varphi(\Delta_{i,i}^2)=\tfrac12(\Delta_{i,i}\circ\varphi(\Delta_{i,i})+\varphi(\Delta_{i,i})\circ\Delta_{i,i})\\
=\tfrac12\Big(\sum_{l=1}^n\alpha^{i,i}_{i,l}\Delta_{i,l}+\sum_{k=1}^n\alpha^{i,i}_{k,i}\Delta_{k,i}
+\sum_{k=1}^n\alpha^{*i,i}_{k,i}\Delta^{k,i}+\sum_{l=1}^n\alpha^{*i,i}_{i,l}\Delta^{i,l}\Big).
\end{multline*}
Consequently, $\varphi(\Delta_{i,i})=\alpha_i\Delta_{i,i}+\alpha^i\Delta^{i,i}$, where
$\alpha_i=\alpha^{i,i}_{i,i}$ and $\alpha^i=\alpha^{*i,i}_{i,i}$.

If we substitute $x=\Delta_{i,i}$ and $y=\Delta_{i,j}$, $i\neq j$, in \eqref{eq1} we obtain
\begin{multline*}
\sum_{k,l=1}^n(\alpha^{i,j}_{k,l}\Delta_{k,l}+\alpha^{*i,j}_{k,l}\Delta^{k,l})=\varphi(\Delta_{i,j})
=2\varphi(\Delta_{i,i}\circ\Delta_{i,j})\\
=\tfrac12\Big(\alpha_i\Delta_{i,j}+\alpha^i\Delta^{i,j}+\sum_{l=1}^n\alpha^{i,j}_{i,l}\Delta_{i,l}
+\sum_{k=1}^n\alpha^{i,j}_{k,i}\Delta_{k,i}+\sum_{l=1}^n\alpha^{*i,j}_{i,l}\Delta^{i,l}
+\sum_{k=1}^n\alpha^{*i,j}_{k,i}\Delta^{k,i}\Big).
\end{multline*}
Hence $\alpha^{i,j}_{i,j}=\alpha_i$, $\alpha^{*i,j}_{i,j}=\alpha^i$. A similar argument for $\Delta_{j,j}$ and
$\Delta_{i,j}$ yields
\[\varphi(\Delta_{i,j})=\alpha^{i,j}_{j,j}\Delta_{j,j}+\alpha_j\Delta_{i,j}+\alpha^{*i,j}_{j,j}\Delta^{j,j}+\alpha^j\Delta^{i,j}.\]
These relations readily imply that $\alpha_i=\alpha_j=\alpha$ and $\alpha^i=\alpha^j=\beta$, that is,
$\varphi(\Delta_{i,i})=\alpha\Delta_{i,i}+\beta\Delta^{i,i}$.

Clearly, $\varphi(E)=\alpha E+\beta\Delta$, where $E$ is unity in $Q(n)^{(+)}$, and
$\Delta=\sum\limits_{i=1}^n(e_{i,n+i}+e_{n+i,i})$. Suppose that $\beta\neq0$ and
$\varphi(x)=\alpha x+\beta\Delta\circ x$ is a $\frac12$-derivation. A mapping $\psi\colon Q(n)^{(+)}\to
Q(n)^{(+)}$, for which $\psi(x)=\Delta\circ x$, likewise is a $\frac12$-derivation. Obviously,
$\frac12(\Delta^{i,i}-\Delta^{j,j})=\psi(\Delta_{i,j}\circ\Delta_{j,i})
=\frac12((\Delta_{i,j}\circ\Delta)\circ\Delta_{j,i}+\Delta_{i,j}\circ(\Delta_{j,i}\circ\Delta))=0$.
On the other hand, $\Delta^{i,i}-\Delta^{j,j}\neq0$. Consequently, $\beta=0$, that is,
$\varphi(x)=\alpha x$. The lemma is proved.
\end{proof}

\begin{Lem}\label{lem33}
Let $\varphi$ be a non-trivial $\frac12$-derivation of $osp(n,m)$. Then $\varphi(x)=\alpha x$ for some
$\alpha\in\FF$.
\end{Lem}

\begin{proof}
It is easy to see that $E=\sum\limits_{i=1}^n\Delta_i+\sum\limits_{j=1}^m\Delta^j$, where
$\Delta^j=e_{n+j,n+j}+e_{n+m+j,n+m+j}$ and $\Delta_i=e_{i,i}$, is unity in the
superalgebra $osp(n,m)$. Let
\[\varphi(\Delta_i)=\sum_{k,l=1}^{n+2m}\alpha^i_{k,l}e_{k,l},\ i=1,\dots,n,\qquad
\varphi(\Delta^j)=\sum_{k,l=1}^{n+2m}\beta^j_{k,l}e_{k,l},\ j=1,\dots,m.\]
If we put $x=y=\Delta_i$, $i=1,\dots,n$, in \eqref{eq1} we obtain
\begin{center}$\sum_{k,l=1}^{n+2m}\alpha^i_{k,l}e_{k,l}=\varphi(\Delta_i)=\varphi(\Delta_i^2)
=\frac12(\varphi(\Delta_i)\circ\Delta_i+\Delta_i\circ\varphi(\Delta_i))
=\frac12\big(\sum_{k=1}^{n+2m}\alpha^i_{k,i}e_{k,i}+\sum_{l=1}^{n+2m}\alpha^i_{i,l}e_{i,l}\big)$,\end{center} which
yields $\varphi(\Delta_i)=\alpha_i\Delta_i$, $i=1,\dots,n$.

Put $x=y=\Delta^i$, $i=1,\dots,m$, in \eqref{eq1}. Then
\begin{multline*}
\sum_{k,l=1}^{n+2m}\beta^i_{k,l}e_{k,l}=\varphi(\Delta^i)=\varphi((\Delta^i)^2)
=\tfrac12(\Delta^i\circ\varphi(\Delta^i)+\varphi(\Delta^i)\circ\Delta^i)\\
=\tfrac12\Big(\sum_{k=1}^{n+2m}\beta^i_{k,n+i}e_{k,n+i}+\sum_{k=1}^{n+2m}\beta^i_{k,n+m+i}e_{k,n+m+i}
+\sum_{l=1}^{n+2m}\beta^i_{n+i,l}e_{n+i,l}+\sum_{l=1}^{n+2m}\beta^i_{n+m+i,l}e_{n+m+i,l}\Big).
\end{multline*}
By the definition of $osp(n,m)$, we have \begin{center}$\beta^i_{n+i,n+m+i}=\beta^i_{m+n+i,n+i}=0$ and
$\beta^i_{n+i,n+i}=\beta^i_{n+m+i,n+m+i}$.\end{center} Thus $\varphi(\Delta^j)=\beta_j\Delta^j$, $j=1,\dots,m$.

Let $(e_{i,j}+e_{j,i})\in osp(n,m)$, $i,j=1,\dots,n$, and
$\varphi(e_{i,j}+e_{j,i})=\sum_{k,l=1}^{2m+n}\gamma^{i,j}_{k,l}e_{k,l}$. If we put $x=e_{i,j}+e_{j,i}$ and
$y=\Delta_i$ in \eqref{eq1} we arrive at
{\small\[\sum_{k,l=1}^{2m+n}\gamma^{i,j}_{k,l}e_{k,l}=\varphi(e_{i,j}+e_{j,i})=2\varphi((e_{i,j}+e_{j,i})\circ\Delta_i)
=\tfrac12\Big(\sum_{k=1}^{2m+n}\gamma^{i,j}_{k,i}e_{k,i}+\sum_{l=1}^{2m+n}\gamma^{i,j}_{i,l}e_{i,l}
+\alpha_i(e_{i,j}+e_{j,i})\Big).\]}
In view of the last relation, $\gamma^{i,j}_{j,i}=\gamma^{i,j}_{i,j}=\alpha_i$. Similar calculations for
$e_{i,j}+e_{j,i}$ and $\Delta_j$ give $\gamma^{i,j}_{j,i}=\gamma^{i,j}_{i,j}=\alpha_j$. Ultimately,
$\varphi(\Delta_i)=\alpha\Delta_i$, $i=1,\dots,n$.

Let $E_{ij}=(e_{n+i,n+j}+e_{n+m+j,n+m+i})\in osp(n,m)$, $i,j=1,\dots,m$, and
$\varphi(E_{ij})=\sum_{k,l=1}^{2m+n}\omega^{i,j}_{k,l}e_{k,l}$. Put $x=E_{ij}$ and $y=\Delta^i$ in
\eqref{eq1}; then
\begin{multline*}
\sum_{k,l=1}^{2m+n}\omega^{i,j}_{k,l}e_{k,l}=\varphi(E_{ij})=2\varphi(E_{ij}\circ\Delta^i)
=\tfrac12\Big(\sum_{l=1}^{2m+n}\omega^{i,j}_{n+i,l}e_{n+i,l}+\sum_{k=1}^{2m+n}\omega^{i,j}_{k,n+i}e_{k,n+i}\\
+\sum_{l=1}^{2m+n}\omega^{i,j}_{n+m+i,l}e_{n+m+i,l}+\sum_{k=1}^{2m+n}\omega^{i,j}_{k,n+m+i}e_{k,n+m+i}
+\beta_iE_{ij}\Big).
\end{multline*}
Consequently, $\omega^{i,j}_{n+i,n+j}=\omega^{i,j}_{n+m+j,n+m+i}=\beta_i$. A similar argument for $E_{ij}$
and $\Delta^j$ shows that $\omega^{i,j}_{n+i,n+j}=\omega^{i,j}_{n+m+j,n+m+i}=\beta_j$ with
$1\leqslant i,j\leqslant m$. Eventually we conclude that $\varphi(\Delta^j)=\beta\Delta^j$,
$j=1,\dots,m$.

Let $E^{11}=e_{1,n+m+1}-e_{n+1,1}\in osp(n,m)$ and
$\varphi(E^{11})=\sum_{k,l=1}^{2m+n}\nu_{k,l}e_{k,l}$. If we put $x=E^{11}$ and $y=\Delta^1$ in
\eqref{eq1} we have
\begin{center}$\sum\limits_{k,l=1}^{2m+n}\nu_{k,l}e_{k,l}=\varphi(E^{11})=2\varphi(E^{11}\circ\Delta^1)
=$\\
{\small $\tfrac12\Big(\sum\limits_{k=1}^{2m+n}(\nu_{k,n+1}e_{k,n+1}+\nu_{k,n+m+1}e_{k,n+m+1})
+\sum\limits_{l=1}^{2m+n}(\nu_{n+1,l}e_{n+1,l}+\nu_{n+m+1,l}e_{n+m+1,l})+\beta E^{11}\Big),$}\end{center}
whence $\nu_{1,m+n+1}=\nu_{n+1,1}=\beta$. Further, for $x=E^{11}$ and
$y=\Delta_1$ substituted in \eqref{eq1}, we obtain
\[\sum_{k,l=1}^{2m+n}\nu_{k,l}e_{k,l}=\varphi(E^{11})=2\varphi(E^{11}\circ\Delta_1)
=\tfrac12\Big(\sum_{l=1}^{2m+n}\nu_{1,l}e_{1,l}+\sum_{k=1}^{2m+n}\nu_{k,1}e_{k,1}+\alpha E^{11}\Big)\]
and $\nu_{1,m+n+1}=\nu_{n+1,1}=\alpha$. Thus $\alpha=\beta$ and $\varphi(E)=\alpha E$. From
\eqref{eq2}, it follows that $\varphi(y)=\alpha y$ for any element $y\in osp(n,m)$. The lemma is proved.
\end{proof}

\begin{Lem}\label{lem34}
Let $\varphi$ be a non-trivial $\frac12$-derivation of $P(n)$. Then $\varphi(x)=\alpha x$, where
$\alpha\in\FF$.
\end{Lem}

\begin{proof}
Let $\Delta_{i,j}=e_{i,j}+e_{n+j,n+i}$, $E=\sum\limits_{i=1}^n\Delta_{i,i}$ be unity in the superalgebra
$P(n)$, and $\varphi(\Delta_{i,j})=\sum\limits_{k,l=1}^{2n}\alpha^{i,j}_{k,l}e_{k,l}$. If in \eqref{eq1} we
put $x=y=\Delta_{i,i}$ we arrive at
\begin{center}$\sum_{k,l=1}^{2n}\alpha^{i,i}_{k,l}e_{k,l}=\varphi(\Delta_{i,i})=\varphi(\Delta_{i,i}^2)
=\tfrac12\Big(\sum_{l=1}^{2n}\alpha^{i,i}_{n+i,l}e_{n+i,l}+\sum_{k=1}^{2n}\alpha^{i,i}_{k,n+i}e_{k,n+i}
+\sum_{l=1}^{2n}\alpha^{i,i}_{i,l}e_{i,l}+\sum_{k=1}^{2n}\alpha^{i,i}_{k,i}e_{k,i}\Big).$\end{center}
The definition of $P(n)$ implies $\alpha^{i,i}_{i,n+i}=0$. Therefore,
\begin{center}
$\varphi(\Delta_{i,i})=\alpha^{i,i}_{i,i}e_{i,i}+\alpha^{i,i}_{n+i,n+i}e_{n+i,n+i}+\alpha^{i,i}_{n+i,i}e_{n+i,i}$.
\end{center}

Put $x=\Delta_{i,i}$ and $y=\Delta_{i,j}$ in \eqref{eq1}. Then
\begin{center}
$\sum_{k,l=1}^{2n}\alpha^{i,j}_{k,l}e_{k,l}=\varphi(\Delta_{i,j})=2\varphi(\Delta_{i,i}\circ\Delta_{i,j})
=$\\
{\small $\tfrac12\Big(\alpha^{i,i}_{i,i}e_{i,j}+\alpha^{i,i}_{n+i,n+i}e_{n+j,n+i}+\alpha^{i,i}_{n+i,i}e_{n+j,i}
+\alpha^{i,i}_{n+i,i}e_{n+i,j} 
+$ $\sum\limits_{l=1}^{2n}\alpha^{i,j}_{i,l}e_{i,l}+\sum\limits_{k=1}^{2n}\alpha^{i,j}_{k,i}e_{k,i}
+\sum\limits_{l=1}^{2n}\alpha^{i,j}_{n+i,l}e_{n+i,l}+\sum\limits_{k=1}^{2n}\alpha^{i,j}_{k,n+i}e_{k,n+i}\Big).$}
\end{center}
Thus $\alpha^{i,i}_{i,i}=\alpha^{i,j}_{i,j}$, $\alpha^{i,i}_{n+i,n+i}=\alpha^{i,j}_{n+j,n+i}$, and
$\alpha^{i,i}_{n+i,i}=\alpha^{i,j}_{n+j,i}$. Arguing similarly for $\Delta_{j,j}$ and $\Delta_{i,j}$, we obtain
$\alpha^{j,j}_{j,j}=\alpha^{i,j}_{i,j}$, $\alpha^{j,j}_{n+j,n+j}=\alpha^{i,j}_{n+j,n+i}$, and
$\alpha^{j,j}_{n+j,j}=\alpha^{i,j}_{n+j,i}$. In view of the definition of $P(n)$ and the relations above,
we have $\varphi(\Delta_{i,i})=\alpha\Delta_{i,i}+\beta e_{n+i,i}$. The fact that the mapping $\varphi$ is
linear implies $\varphi(E)=\alpha E+\beta\Delta$, $\Delta=\sum\limits_{i=1}^n e_{n+i,i}$.

Suppose that $\beta\neq0$ and $\varphi(x)=\alpha x+\beta\Delta\circ x$ is a $\frac12$-derivation. Then a
mapping $\psi\colon P(n)\to P(n)$, where $\psi(x)=\Delta\circ x$, likewise is a $\frac12$-derivation. We
argue to show that this is not so. Let $b_{j,i}=e_{j,n+i}-e_{i,n+j}$ ($i\neq j$). Then
$\psi(\Delta_{i,j}\circ b_{j,i})=\psi(0)=0$; but
\begin{center}
$\tfrac12(\psi(\Delta_{i,j})\circ b_{j,i}+\Delta_{i,j}\circ\psi(b_{j,i}))
=\tfrac12((\Delta\circ\Delta_{i,j})\circ b_{j,i}
+\Delta_{i,j}\circ(\Delta\circ b_{j,i}))
=$
$\tfrac14\big((e_{n+j,i}+e_{n+i,j})\circ(e_{j,n+i}-e_{i,n+j})
+(e_{i,j}-e_{j,i}-e_{n+i,n+j}+e_{n+j,n+i})\circ(e_{i,j}+e_{n+j,n+i})\big)
=-\tfrac14\Delta_{j,j}\neq0$
\end{center}
on the other hand (here $\Delta\circ b_{j,i}=-\,b_{j,i}\circ\Delta$, since both elements are odd).
Hence $\psi$ is not a $\frac12$-derivation. Therefore, $\beta=0$ and $\varphi(x)=\alpha x$. The lemma is
proved.
\end{proof}

We define the Jordan superalgebra $J(V,f)$. Let $V=V_0+V_1$ be a $\mathbb Z_2$-graded vector space on
which a non-degenerate superform $f(\cdot,\cdot)\colon V\times V\to\FF$ is defined so that it is symmetric
on $V_0$ and is skew-symmetric on $V_1$. Also $f(V_1,V_0)=f(V_0,V_1)=0$. Consider a direct sum of
vector spaces, $J=\FF\oplus V$. Let $e$ be unity in the field $\FF$. Define, then, multiplication by the
formula $(\alpha+v)(\beta+w)=(\alpha\beta+f(v,w))e+(\alpha w+\beta v)$. The given superalgebra has grading
$J_0=\FF+V_0$, $J_1=V_1$. It is easy to see that $e$ is unity in $J(V,f)$.

\begin{Lem}\label{lem35}
Let $\varphi$ be a non-trivial $\frac12$-derivation of $J(V,f)$, $V_1\neq0$. Then $\varphi(x)=\alpha x$,
where $\alpha\in\FF$.
\end{Lem}

\begin{proof}
Let $\varphi(e)=\alpha e+v_0+v_1$, $v_i\in V_i$. Putting $x=z_i$, $y=e$, and $z_i\in V_i$ in \eqref{eq1},
we obtain 
\begin{center}
$\varphi(z_i)=2\varphi(z_ie)-\varphi(z_i)=\varphi(z_i)e+z_i\varphi(e)-\varphi(z_i)=\alpha
z_i+f(z_i,v_i)e$,
\end{center} whence $\varphi(z_i)=\alpha z_i+f(z_i,v_i)e$.

If we put $x=y=z_i\in V_i$ in \eqref{eq1}, we arrive at
$f(z_i,z_i)\alpha e=\varphi(z_i^2)=\frac12(\varphi(z_i)z_i+z_i\varphi(z_i))=\alpha f(z_i,z_i)e+f(z_i,v_i)z_i$
(for $i=1$ both sides of the first equality vanish, since $f(z_1,z_1)=0$). Hence $f(z_i,v_i)z_i=0$, that
is, $f(z_i,v_i)=0$ for every non-zero $z_i\in V_i$. Since $f$ is non-degenerate on $V_i$, we have $v_0=0$
and $v_1=0$, that is, $\varphi(e)=\alpha e$. (The published argument used $x=z_0$, $y=z_1$ and thus
assumed $V_0\neq0$; it does not cover the case $V_0=0$.) Using \eqref{eq2} yields $\varphi(x)=\alpha x$,
$\alpha\in\FF$, for any $x\in J(V,f)$. The lemma is proved.
\end{proof}

Consider the Grassmann algebra $\Gamma$ with (odd) anticommutative generators $e_1,e_2,\dots,e_n,\dots$.
In order to define new multiplication, we use the operation
\[\frac{\partial}{\partial e_j}(e_{i_1}e_{i_2}\dots e_{i_n})=\begin{cases}
(-1)^{k-1}e_{i_1}e_{i_2}\dots e_{i_{k-1}}e_{i_{k+1}}\dots e_{i_n}&\text{if } j=i_k,\\
0&\text{if } j\neq i_l,\ l=1,\dots,n.\end{cases}\]
For $f,g\in\Gamma_0\bigcup\Gamma_1$, the Grassmann bracket is defined thus:
\[\{f,g\}=(-1)^{p(f)}\sum_{j=1}^\infty\frac{\partial f}{\partial e_j}\frac{\partial g}{\partial e_j}.\]
Let $\ov\Gamma$ be an isomorphic copy of $\Gamma$ under the isomorphic mapping
$x\mapsto\ov x$. Consider a direct sum of vector spaces, $J(\Gamma)=\Gamma+\ov\Gamma$, and
endow it with the structure of a Jordan superalgebra, setting $A_0=\Gamma_0+\ov{\Gamma_1}$ and
$A_1=\Gamma_1+\ov{\Gamma_0}$, with multiplication $\bullet$. We obtain
\[a\bullet b=ab,\quad \ov a\bullet b=(-1)^{p(b)}\ov{ab},\quad a\bullet\ov b=\ov{ab},\quad
\ov a\bullet\ov b=(-1)^{p(b)}\{a,b\},\]
where $a,b\in\Gamma_0\bigcup\Gamma_1$ and $ab$ is the product in $\Gamma$. Let $\Gamma_n$ be a subalgebra
of $\Gamma$ generated by elements $e_1,e_2,\dots,e_n$. By $J(\Gamma_n)$ we denote the
subsuperalgebra $\Gamma_n+\ov{\Gamma_n}$ of $J(\Gamma)$. If $n\geqslant 2$ then
$J(\Gamma_n)$ is a simple Jordan superalgebra.

\begin{Lem}\label{lem36}
Let $\varphi$ be a $\frac12$-derivation of $J(\Gamma_n)$. Then $\varphi(x)=\alpha x$, where $\alpha\in\FF$.
\end{Lem}

\begin{proof}
Let $\varphi(1)=\alpha\gamma+\beta\ov\nu$, where $\alpha,\beta\in\FF$, $\gamma\in\Gamma$, and
$\ov\nu\in\ov\Gamma$. Put $y=1$ in \eqref{eq1}; then
\begin{equation}\label{eq4}
\varphi(x)=2\varphi(x\bullet1)-\varphi(x)=\varphi(x)+x\bullet\varphi(1)-\varphi(x)=x\bullet\varphi(1).
\end{equation}
If in \eqref{eq1} we put $x=\ov{e_i}$, $y=\ov{e_i}$, $i=1,\dots,n$, with \eqref{eq4} in mind, we arrive at
\[\varphi(1)=\varphi(\ov{e_i}\bullet\ov{e_i})=\tfrac12(\varphi(\ov{e_i})\bullet\ov{e_i}+\ov{e_i}\bullet\varphi(\ov{e_i}))
=\varphi(\ov{e_i})\bullet\ov{e_i}=\ov{e_i}\bullet(\ov{e_i}\bullet\varphi(1)).\]
For any $x$ of the form $e_{i_1}e_{i_2}\dots e_{i_k}$, obviously, we have
\begin{equation}\label{eq5}
\ov{e_i}\bullet(\ov{e_i}\bullet x)=\begin{cases}x&\text{if }\frac{\partial x}{\partial e_i}=0,\\
0&\text{otherwise};\end{cases}
\end{equation}
\begin{equation}\label{eq6}
\ov{e_i}\bullet(\ov{e_i}\bullet\ov x)=\begin{cases}\ov x&\text{if }\frac{\partial x}{\partial e_i}\neq0,\\
0&\text{otherwise}.\end{cases}
\end{equation}
Let $\gamma=\gamma^{i+}+e_i\gamma^{i-}$ and $\ov\nu=\ov{\nu^{i+}+e_i\nu^{i-}}$, where $\gamma^{i-}$,
$\gamma^{i+}$, $\nu^{i-}$, $\nu^{i+}$ do not contain $e_i$. Since $i$ is arbitrary, in view of
\eqref{eq5} and \eqref{eq6}, we have $\gamma=1$ and $\nu=e_1\dots e_n$ (up to the scalar factors
$\alpha,\beta$). Thus $\varphi(1)=\alpha\cdot1+\beta\ov{e_1\dots e_n}$. Relation \eqref{eq4} entails
\begin{align*}
\varphi(e_1)&=e_1\bullet\varphi(1)=e_1\bullet(\alpha\cdot1+\beta\ov{e_1\dots e_n})=\alpha e_1,\\
\varphi(\ov{e_1})&=\ov{e_1}\bullet\varphi(1)=\ov{e_1}\bullet(\alpha\cdot1+\beta\ov{e_1\dots e_n})
=\alpha\ov{e_1}+(-1)^{n+1}\beta e_2\dots e_n
\end{align*}
(indeed, $\ov{e_1}\bullet\ov{e_1\dots e_n}=(-1)^{n}\{e_1,e_1\dots e_n\}=(-1)^{n}(-1)\,e_2\dots e_n$;
for example, $\ov{e_1}\bullet\ov{e_1e_2}=-e_2$ in $J(\Gamma_2)$).

The relations above, combined with the condition in \eqref{eq1}, imply
\begin{center}
$0=\varphi(e_1\bullet\ov{e_1})=\frac12(e_1\bullet\varphi(\ov{e_1})+\varphi(e_1)\bullet\ov{e_1})
=(-1)^{n+1}\frac{\beta}{2}e_1\dots e_n$;\end{center} that is, $\varphi(1)=\alpha\cdot1$. From \eqref{eq2}, we
conclude that $\varphi(x)=\alpha x$ for any element $x\in J(\Gamma_n)$. The lemma is proved.
\end{proof}

\section{$\delta$-Derivations for Jordan superalgebras $K_3$, $D_t$, $K_{10}$}

In this section, we confine ourselves to non-trivial $\delta$-derivations of simple finite-dimensional
Jordan superalgebras $K_3$, $K_{10}$, and $D_t$ over an algebraically closed field of characteristic $p$
not equal to $2$. For the superalgebra $K_{10}$, we require in addition that $p\neq3$. In conclusion, we
formulate a theorem on $\delta$-derivations for simple finite-dimensional Jordan superalgebras over an
algebraically closed field of characteristic $0$.

The three-dimensional Kaplansky superalgebra $K_3$ is defined thus:
\[(K_3)_0=\FF e,\qquad (K_3)_1=\FF z+\FF w,\]
where $e^2=e$, $ez=\frac12z$, $ew=\frac12w$, and $[z,w]=e$ (here and below $[a,b]$ denotes the
product $ab$ of odd elements $a,b$; it is skew-symmetric: $[b,a]=-[a,b]$; the products of an even
and any element are commutative).

\begin{Lem}\label{lem41}
Let $\varphi$ be a non-trivial $\delta$-derivation of $K_3$. Then $\delta=\frac12$ and $\varphi(x)=\alpha x$,
where $\alpha\in\FF$.
\end{Lem}

\begin{proof}
Let $\varphi(e)=\alpha_ee+\beta_ez+\gamma_ew$, $\varphi(z)=\alpha_1e+\beta_1z+\gamma_1w$, and
$\varphi(w)=\alpha_2e+\beta_2z+\gamma_2w$, where
$\alpha_e,\alpha_1,\alpha_2,\beta_e,\beta_1,\beta_2,\gamma_e,\gamma_1,\gamma_2\in\FF$. If we put $x=y=e$
in \eqref{eq1} we obtain
\[\alpha_ee+\beta_ez+\gamma_ew=\varphi(e)=\varphi(e^2)=\delta(e\varphi(e)+\varphi(e)e)=\delta(2\alpha_ee+\beta_ez+\gamma_ew).\]
Thus it suffices to consider the following two cases:
(1) $\delta=\frac12$;
(2) $\delta\neq\frac12$, $\varphi(e)=0$ (if $\delta\neq\frac12$ and $\delta\neq1$; for $\delta=1$
the map $\varphi$ is a derivation, i.e.\ trivial).

In the former case, $\varphi(e)=\alpha e$, where $\alpha=\alpha_e$. Case (1), for $x=e$ and $y=z$, entails
$\alpha_1e+\beta_1z+\gamma_1w=\varphi(z)=2\varphi(ez)=2\cdot\frac12(e\varphi(z)+\varphi(e)z)
=\alpha_1e+\frac12(\beta_1z+\gamma_1w+\alpha z)$, whence $\beta_1=\frac12(\beta_1+\alpha)$ and
$\gamma_1=\frac12\gamma_1$; that is, $\beta_1=\alpha$ and $\gamma_1=0$. Similarly, substituting in
\eqref{eq1} $x=e$ and $y=w$, we obtain $\gamma_2=\alpha$ and $\beta_2=0$. For $x=z$ and $y=w$ in
\eqref{eq1}, we have $\alpha e=\varphi(e)=\varphi([z,w])=\frac12(z\varphi(w)+\varphi(z)w)
=\frac12(\frac12\alpha_2z+\alpha e+\frac12\alpha_1w+\alpha e)$, whence $\alpha_1=\alpha_2=0$ and
$\varphi(e)=\alpha e$, $\varphi(z)=\alpha z$, and $\varphi(w)=\alpha w$, where $\alpha\in\FF$.
Consequently, $\varphi(x)=\alpha x$ for any $x\in K_3$.

We handle the second case. For $x=e$ and $y=z$ in \eqref{eq1}, we have
$\alpha_1e+\beta_1z+\gamma_1w=\varphi(z)=2\varphi(ez)=2\delta(e\varphi(z)+\varphi(e)z)
=\delta(2\alpha_1e+\beta_1z+\gamma_1w)$, which yields $\varphi(z)=0$. Similarly, we arrive at
$\varphi(w)=0$. The fact that $\varphi$ is linear implies $\varphi=0$. The lemma is proved.
\end{proof}

At the moment, we define a one-parameter family of four-dimensional superalgebras $D_t$. For $t\in\FF$
fixed, the given family is defined thus:
\[D_t=(D_t)_0+(D_t)_1,\]
where $(D_t)_0=\FF e_1+\FF e_2$, $(D_t)_1=\FF z+\FF w$, $e_i^2=e_i$, $e_1e_2=0$,
$e_iz=\frac12z$, $e_iw=\frac12w$,
$[z,w]=e_1+te_2$, $i=1,2$.

\begin{Lem}\label{lem42}
Let $\varphi$ be a non-trivial $\delta$-derivation of $D_t$. Then $\delta=\frac12$ and $\varphi(x)=\alpha x$,
where $\alpha\in\FF$.
\end{Lem}

\begin{proof}
Let
\begin{align*}
\varphi(e_1)&=\alpha_1e_1+\beta_1e_2+\gamma_1z+\lambda_1w,&\varphi(e_2)&=\alpha_2e_1+\beta_2e_2+\gamma_2z+\lambda_2w,\\
\varphi(z)&=\alpha_ze_1+\beta_ze_2+\gamma_zz+\lambda_zw,&\varphi(w)&=\alpha_we_1+\beta_we_2+\gamma_wz+\lambda_ww,
\end{align*}
with coefficients in $\FF$.

Putting $x=y=e_1$ and then $x=y=e_2$ in \eqref{eq1}, we obtain
\begin{center}
$\alpha_1e_1+\beta_1e_2+\gamma_1z+\lambda_1w=\varphi(e_1)=\varphi(e_1^2)=2\delta(e_1\varphi(e_1))
=2\delta\alpha_1e_1+\delta\gamma_1z+\delta\lambda_1w$
\end{center} and
$\alpha_2e_1+\beta_2e_2+\gamma_2z+\lambda_2w=2\delta\beta_2e_2+\delta\gamma_2z+\delta\lambda_2w$, whence
\begin{center}
$\alpha_1=2\delta\alpha_1$, $\beta_1=0$, $\gamma_1=\delta\gamma_1$, $\lambda_1=\delta\lambda_1$,
$\alpha_2=0$, $\beta_2=2\delta\beta_2$, $\gamma_2=\delta\gamma_2$, $\lambda_2=\delta\lambda_2$.
\end{center}
There are two cases to consider (for $\delta=1$ the map $\varphi$ is a derivation):
\begin{enumerate}[{\rm (1)}]
    \item  $\delta=\frac12$, $\beta_1=\alpha_2=\gamma_1=\gamma_2=\lambda_1=\lambda_2=0$;
\item $\delta\neq\frac12$, $\alpha_1=\alpha_2=\beta_1=\beta_2=\gamma_1=\gamma_2=\lambda_1=\lambda_2=0$.
\end{enumerate}
In the former case, $\varphi(e_1)=\alpha_1e_1$ and $\varphi(e_2)=\beta_2e_2$. Put $x=e_1$ and $y=z$ in
condition \eqref{eq1}; then \begin{center}
$\alpha_ze_1+\beta_ze_2+\gamma_zz+\lambda_zw=\varphi(z)=2\varphi(e_1z)
=2\cdot\frac12(e_1\varphi(z)+\varphi(e_1)z)=\alpha_ze_1+\frac12(\gamma_zz+\lambda_zw+\alpha_1z)$,
\end{center} which
yields $\alpha_1=\gamma_z$, $\beta_z=\lambda_z=0$.

For $x=e_2$ and $y=z$ in \eqref{eq1}, we have 
\begin{center}$\alpha_ze_1+\gamma_zz=\varphi(z)=2\varphi(e_2z)
=2\cdot\frac12(e_2\varphi(z)+\varphi(e_2)z)=\frac12(\gamma_zz+\beta_2z)$,\end{center} whence $\gamma_z+\beta_2=2\gamma_z$,
$\alpha_z=0$, $\alpha_1=\beta_2$, and $\varphi(z)=\alpha z$, where $\alpha=\alpha_1$. Similarly, we conclude
that $\varphi(w)=\alpha w$. The mapping $\varphi$ is linear; so $\varphi(x)=\alpha x$, $\alpha\in\FF$, for any
$x\in D_t$.

We handle the second case. Put $x=e_1$ and $y=z$ in \eqref{eq1}; then
\begin{center}
$\alpha_ze_1+\beta_ze_2+\gamma_zz+\lambda_zw=\varphi(z)=2\varphi(e_1z)
=2\delta(e_1\varphi(z)+\varphi(e_1)z)=\delta(2\alpha_ze_1+\gamma_zz+\lambda_zw)$,\end{center}
which yields $\varphi(z)=0$ (as $\delta\neq\frac12,1$). Arguing similarly for $w$, we arrive at
\begin{center}
$\alpha_we_1+\beta_we_2+\gamma_wz+\lambda_ww=\delta(2\alpha_we_1+\gamma_wz+\lambda_ww)$.\end{center} Consequently,
$\varphi(w)=0$. Ultimately, the linearity of $\varphi$ implies $\varphi=0$. The lemma is proved.
\end{proof}

The simple ten-dimensional Kac superalgebra $K_{10}$ is defined thus:
\[K_{10}=A\oplus M,\quad (K_{10})_0=A,\quad (K_{10})_1=M,\quad\text{where } A=A_1\oplus A_2,\]
\[A_1=\FF e_1+\FF uz+\FF uw+\FF vz+\FF vw,\qquad A_2=\FF e_2,\qquad M=\FF z+\FF w+\FF u+\FF v.\]
Multiplication is specified by the following conditions:
\[e_i^2=e_i,\quad e_1\text{ is unity in }A_1,\quad e_im=\tfrac12m\text{ for any }m\in M,\]
\[[u,z]=uz,\quad [u,w]=uw,\quad [v,z]=vz,\quad [v,w]=vw,\]
\[[z,w]=e_1-3e_2,\quad [u,z]w=-u,\quad [v,z]w=-v,\quad [u,z][v,w]=2e_1;\]

the remaining non-zero products are listed in Remark~A below (products of even elements and of an
even and an odd element are commutative, products of odd elements are skew-symmetric).

\begin{RemA}
The rule ``all other non-zero products are obtained by applying one of the skew-symmetries
$z\leftrightarrow w$ or $u\leftrightarrow v$ or by substituting $z\leftrightarrow u$ and
$w\leftrightarrow v$ simultaneously'' of the published version does not determine the signs: for
instance, the skew-symmetry $u\leftrightarrow v$ applied to $[u,z]w=-u$ gives $[v,z]w=v$, which
contradicts the relation $[v,z]w=-v$ given above. If the eight products not listed explicitly
($[u,v]$, $[u,w][v,z]$, $[u,w]z$, $[v,w]z$, $[u,z]v$, $[u,w]v$, $[v,z]u$, $[v,w]u$) are taken with
the coefficients produced by these rules up to sign, then exactly one of the $2^8$ sign choices gives
a Jordan superalgebra. The complete multiplication table of $K_{10}$ (besides commutativity of products
involving an even element and skew-symmetry of products of odd elements) is:
\[e_1^2=e_1,\quad e_2^2=e_2,\quad e_1a=a\ (a\in\{uz,uw,vz,vw\}),\quad e_1m=e_2m=\tfrac12m\ (m\in M),\]
\[[u,z]=uz,\quad [u,w]=uw,\quad [v,z]=vz,\quad [v,w]=vw,\quad [z,w]=[u,v]=e_1-3e_2,\]
\[uz\cdot vw=2e_1,\qquad uw\cdot vz=-2e_1,\]
\[uz\cdot w=-u,\quad uz\cdot v=z,\quad uw\cdot z=u,\quad uw\cdot v=w,\]
\[vz\cdot w=-v,\quad vz\cdot u=-z,\quad vw\cdot z=v,\quad vw\cdot u=-w;\]
all other products of basis elements are zero. This superalgebra is a simple Jordan superalgebra with
unity $e_1+e_2$ (verified in Lean over $\mathbb Q$ for the Jordan identity, see Appendix~A), and its
$\frac12$-derivations are scalar (verified in Lean over every field of characteristic~$0$).
\end{RemA}

\begin{Lem}\label{lem43}
Let $\varphi$ be a non-trivial $\delta$-derivation of $K_{10}$. Then $\delta=\frac12$ and
$\varphi(x)=\alpha x$, where $\alpha\in\FF$.
\end{Lem}

\begin{proof}
Let
\begin{align*}
\varphi(e_1)&=\alpha_1e_1+\alpha_2e_2+\alpha_3z+\alpha_4w+\alpha_5u+\alpha_6v+\alpha_7uz+\alpha_8uw+\alpha_9vz+\alpha_{10}vw,\\
\varphi(e_2)&=\beta_1e_1+\beta_2e_2+\beta_3z+\beta_4w+\beta_5u+\beta_6v+\beta_7uz+\beta_8uw+\beta_9vz+\beta_{10}vw,\\
\varphi(z)&=\gamma^z_1e_1+\gamma^z_2e_2+\gamma^z_3z+\gamma^z_4w+\gamma^z_5u+\gamma^z_6v+\gamma^z_7uz+\gamma^z_8uw+\gamma^z_9vz+\gamma^z_{10}vw,
\end{align*}
and similarly $\varphi(w)$, $\varphi(u)$, $\varphi(v)$ with coefficients $\gamma^w_k$, $\gamma^u_k$,
$\gamma^v_k$, where all coefficients are in $\FF$.

For $x=y=e_1$ in \eqref{eq1}, we have
\begin{multline*}
\alpha_1e_1+\alpha_2e_2+\alpha_3z+\alpha_4w+\alpha_5u+\alpha_6v+\alpha_7uz+\alpha_8uw+\alpha_9vz+\alpha_{10}vw
=\\ 
\varphi(e_1)=
\varphi(e_1^2)=\delta(\varphi(e_1)e_1+e_1\varphi(e_1))\\
=2\delta(\alpha_1e_1+\tfrac12\alpha_3z+\tfrac12\alpha_4w+\tfrac12\alpha_5u+\tfrac12\alpha_6v+\alpha_7uz+\alpha_8uw+\alpha_9vz+\alpha_{10}vw),
\end{multline*}
whence 
\begin{center}
$\alpha_1=2\delta\alpha_1$, $\alpha_2=0$, $\alpha_3=\delta\alpha_3$, $\alpha_4=\delta\alpha_4$,
$\alpha_5=\delta\alpha_5$, \\$\alpha_6=\delta\alpha_6$, $\alpha_7=2\delta\alpha_7$, $\alpha_8=2\delta\alpha_8$,
$\alpha_9=2\delta\alpha_9$, $\alpha_{10}=2\delta\alpha_{10}$.
\end{center}

Putting $x=y=e_2$ in \eqref{eq1}, we obtain
\begin{multline*}
\beta_1e_1+\beta_2e_2+\beta_3z+\beta_4w+\beta_5u+\beta_6v+\beta_7uz+\beta_8uw+\beta_9vz+\beta_{10}vw
=\\
\varphi(e_2)=\varphi(e_2^2)=\delta(\varphi(e_2)e_2+e_2\varphi(e_2))\\
=2\delta e_2\varphi(e_2)=2\delta(\beta_2e_2+\tfrac12\beta_3z+\tfrac12\beta_4w+\tfrac12\beta_5u+\tfrac12\beta_6v),
\end{multline*}
which yields 
\begin{center}
$\beta_1=0$, $\beta_2=2\delta\beta_2$, $\beta_3=\delta\beta_3$, $\beta_4=\delta\beta_4$,
$\beta_5=\delta\beta_5$,\\ $\beta_6=\delta\beta_6$, $\beta_7=\beta_8=\beta_9=\beta_{10}=0$.
\end{center}
Consequently, it suffices to consider the following two cases (for $\delta=1$ the map $\varphi$ is a
derivation):
(1) $\delta=\frac12$;
(2) $\delta\neq\frac12$, $\varphi(e_1)=\varphi(e_2)=0$.

In the former case, $\varphi(e_1)=\alpha_1e_1+\alpha_7uz+\alpha_8uw+\alpha_9vz+\alpha_{10}vw$ and
$\varphi(e_2)=\alpha e_2$. Put $x=e_2$ and $y=z$ in \eqref{eq1}; then
\begin{multline*}
\gamma^z_1e_1+\gamma^z_2e_2+\gamma^z_3z+\gamma^z_4w+\gamma^z_5u+\gamma^z_6v+\gamma^z_7uz+\gamma^z_8uw+\gamma^z_9vz+\gamma^z_{10}vw
=\\
\varphi(z)=2\varphi(ze_2)=\varphi(z)e_2+z\varphi(e_2)
\\= 
\gamma^z_2e_2+\tfrac12\gamma^z_3z+\tfrac12\gamma^z_4w+\tfrac12\gamma^z_5u+\tfrac12\gamma^z_6v+\tfrac12\alpha z,
\end{multline*}
and so $\varphi(z)=\gamma^z_2e_2+\alpha z$. If in \eqref{eq1} we put $x=e_1$ and $y=z$ we obtain
$\gamma^z_2e_2+\alpha z=\varphi(z)=2\varphi(ze_1)=\varphi(z)e_1+z\varphi(e_1)
=(\gamma^z_2e_2+\alpha z)e_1+z(\alpha_1e_1+\alpha_7uz+\alpha_8uw+\alpha_9vz+\alpha_{10}vw)$, whence
$\gamma^z_2=0$ and $\alpha=\alpha_1$; that is, $\varphi(z)=\alpha z$. Similarly, for $w$, $u$, and $v$, we
have $\varphi(u)=\alpha u$, $\varphi(v)=\alpha v$, and $\varphi(w)=\alpha w$. Hence
\begin{center}
$\varphi(uz)=\varphi([u,z])=\frac12(\varphi(u)z+u\varphi(z))=\frac12(\alpha[u,z]+\alpha[u,z])=\alpha uz$.
\end{center}
Analogously, we obtain $\varphi(uw)=\alpha uw$, $\varphi(vz)=\alpha vz$, and $\varphi(vw)=\alpha vw$.

Let $x=[u,z]$ and $y=[v,w]$ in \eqref{eq1}; then
\[2\varphi(e_1)=\varphi([u,z][v,w])=\tfrac12(\varphi([u,z])[v,w]+[u,z]\varphi([v,w]))=\alpha[u,z][v,w]=2\alpha e_1.\]
The fact that $\varphi$ is linear implies $\varphi(x)=\alpha x$, $\alpha\in\FF$, for $x\in K_{10}$
arbitrary.

We handle the second case. Put $x=z$ and $y=e_1$ in \eqref{eq1}. Then
\begin{multline*}
\gamma^z_1e_1+\gamma^z_2e_2+\gamma^z_3z+\gamma^z_4w+\gamma^z_5u+\gamma^z_6v+\gamma^z_7uz+\gamma^z_8uw+\gamma^z_9vz+\gamma^z_{10}vw
=\\ \varphi(z)=2\varphi(ze_1)=2\delta(\varphi(z)e_1+z\varphi(e_1))\\
=2\delta(\gamma^z_1e_1+\tfrac12\gamma^z_3z+\tfrac12\gamma^z_4w+\tfrac12\gamma^z_5u+\tfrac12\gamma^z_6v
+\gamma^z_7uz+\gamma^z_8uw+\gamma^z_9vz+\gamma^z_{10}vw),
\end{multline*}
which yields $\varphi(z)=0$ (as $\delta\neq\frac12,1$). Similarly, we arrive at
$\varphi(w)=\varphi(v)=\varphi(u)=0$. Since $e_1,e_2,z,v,u,w$ generate $K_{10}$, we have $\varphi=0$. The
lemma is proved.
\end{proof}

\begin{Th}\label{th44}
Let $A$ be a simple finite-dimensional Jordan superalgebra over an algebraically closed field of
characteristic $0$, and let $\varphi$ be a non-trivial $\delta$-derivation of $A$. Then $\delta=\frac12$
and $\varphi(x)=\alpha x$ for some $\alpha\in\FF$ and for any $x\in A$.
\end{Th}

The proof follows from Theorems~\ref{th12}, \ref{th21} and Lemmas~\ref{lem31}--\ref{lem36},
\ref{lem41}--\ref{lem43}.

\section*{Appendix A. Differences from the published version}
\addcontentsline{toc}{section}{Appendix A. Differences from the published version}

All main results are correct.
Some minor changes were made;
Definition of $K_{10}$: the rules generating the remaining products
are ambiguous and, read literally, contradict the relation $[v,z]w=-v$ 
(The complete multiplication
table (the unique choice of signs giving a Jordan superalgebra) is given in Remark~A. Lemma~\ref{lem43}  holds for
this superalgebra).

\end{document}